\documentclass[12pt]{article}
\usepackage{amssymb,amsmath}
\usepackage{mathrsfs}

\newtheorem{Thm}{Theorem}

\begin{document}

\begin{center}
\large{\bf STRONG AND WEAK CONVERGENCE OF POPULATION SIZE IN SUPERCRITICAL CATALYTIC BRANCHING PROCESS}
\end{center}
\vskip0,5cm
\begin{center}
Ekaterina Vl. Bulinskaya\footnote{ \emph{Email address:} {\tt
bulinskaya@yandex.ru}}$^,$\footnote{The work is partially supported
by Dmitry Zimin Foundation ``Dynasty'' and RFBR grant 14-01-00318.}
\vskip0,2cm \emph{Lomonosov Moscow State University}
\end{center}
\vskip1cm

\begin{abstract}

A general model of catalytic branching process (CBP) with any finite number of catalysis centers in a discrete space is studied.
More exactly, it is assumed that particles move in this space according to a specified Markov chain and they may produce offspring only in the presence of catalysts located at fixed points. The asymptotic (in time) behavior of the total number of particles as well as the local particles numbers is investigated. The problems of finding the global extinction probability and local extinction probability are solved. Necessary and sufficient conditions are established for phase of pure global survival and strong local survival. Under wide conditions the limit theorems for the normalized total and local particles numbers in supercritical CBP are proved in the sense of almost surely convergence as well as with respect to convergence in distribution. Generalizations of a number of previous results are obtained as well. In the proofs the main role is played by recent results by the author devoted to classification of CBP and the moment analysis of the total and local particles numbers in CBP.

\vskip0,5cm {\it Keywords and phrases}: catalytic branching process,
extinction probability, strong local survival, pure global phase, limit theorems,
multi-type Bellman-Harris process.

\vskip0,5cm 2010 {\it AMS classification}: 60J80, 60J27.

\end{abstract}

Theory of branching processes is a classical part of probability
theory (see, e.g., \cite{Sev}). Its origin goes back to the middle
of XIX century when F.Galton and G.W.Watson proposed a model
elucidating the extinction of family names. In the framework of this
model the primary branching process called Galton-Watson process
arose. It was found that the value of  extinction probability as
well as the mean population size depend essentially on whether the
mean offspring number of a population representative is greater
than, equal to or less than 1. Galton-Watson branching process
having that mean offspring number is called supercritical, critical
and subcritical, respectively. Only for supercritical Galton-Watson
process the extinction probability is less than 1 whereas the mean
population size tends to infinity as time grows. Therefore, after
solving the problem of finding the extinction probability the
question arises how fast the population grows in the case of its
survival. The answer to that question is given by limit theorems for
the population size. The results related to the Galton-Watson
process are contained, e.g., in \cite{Vatutin_08}, sections 1-6 and
9.

We are interested in a more involved model than the Galton-Watson
process, namely catalytic branching process (CBP) with an arbitrary
finite number of catalysis centers. Its specifics consists in
possibility for the population representatives (which further will
be called particles) not only to produce offsprings but also to move
in the space. Moreover, we assume that the particles produce
offsprings exclusively in the presence of catalysts which are
located in a finite set of the space points. For CBP it is natural
to raise the question not only on the global particles extinction
but also on the local extinction as well as to investigate the
asymptotic behavior of total and local particles numbers. The
present paper is devoted to these problems. We consider a discrete
space where the particles move, although branching processes in a
continuous space are of doubtless interest as well (see, e.g.,
\cite{Bocharov_Harris_2014}).

The study of general model of CBP with any finite number of
catalysis centers in a discrete space is initiated in
\cite{B_TVP_14}. Particular cases of this model were analyzed in
many publications, e.g., in \cite{ABY_98},
\cite{Carmona_Hu_14}, \cite{DR_13}, \cite{VT_MT_11} and others. In dependence on the
value of the Perron root of a certain matrix a classification of CBP
as supercritical, critical and subcritical processes was proposed in
\cite{B_TVP_14}. The naturalness of the classification was confirmed
by implemented by the author asymptotic analysis of the moments of
the total and local particles numbers existing in the process at
time $t$, as $t\to\infty$. In the present work the problem of
finding the probabilities of global and local population extinction
is solved and new limit theorems in strong (Theorem \ref{strong})
and weak (Theorem \ref{weak}) forms are established for the total
and local particles numbers in a supercritical CBP. The results
obtained generalize a number of previous ones proved, e.g., in
papers  \cite{Carmona_Hu_14} (Lemma 5.1) and \cite{Y_MZ_12}. So,
limit distributions for particles numbers in the model of
supercritical branching random walk on integer lattice
$\mathbb{Z}^d$, $d\in\mathbb{N}$, with a finite number of particles
generation centers are found in \cite{Y_MZ_12}. There not only
existence of all moments of an offsprings number for any particle
was assumed but also a certain rate of growth of such moments in
dependence on its order. In our work, besides considering a more
general model we also impose extremely weak restrictions on the
moments of the offsprings size for each particle and, moreover,
study the limit behavior of the particles numbers in the sense of
almost surely convergence. These results are not only of
self-contained interest but will be also applied in further
investigations of the character of the particles population
propagation in CBP. It is worthwhile to note that the mentioned
achievements in the domain of limit theorems for supercritical CBP
are due to the foundation laid in a recent work \cite{B_TVP_14}.
Firstly, there were introduced auxiliary Bellman-Harris branching
processes which mainly allows us to reduce the study of CBP to
classic results of modern theory of branching processes. Secondly,
the moment analysis of the total and local particles numbers was
implemented in \cite{B_TVP_14} and the proofs of results of the
present paper substantially bear on it.

Let us give a formal description of CBP. At the initial time $t=0$
there is a single particle, its movement on a specified finite or
countable set $S$ is governed by continuous-time Markov chain
$\eta=\{\eta(t),t\geq0\}$ with infinitesimal matrix
$Q=(q(x,y))_{x,y\in S}$. Hitting a finite catalysts set
${W=\{w_1,\ldots,w_N\}\subset S}$, e.g., at point $w_k$, the
particle spends there random time distributed according to the
exponential law with parameter $\beta_k>0$. Then it either produces
offsprings or leaves the element $w_k$ with probabilities $\alpha_k$
and $1-\alpha_k$ ($0\leq\alpha_k<1$), respectively. If the particle
produces offsprings (at $w_k$), it dies instantly turning into a
random number $\xi_{k}$ of offsprings located at the same point
$w_k$. If the particle leaves $w_k$, it jumps to point $y\neq w_k$
with probability $-q(w_k,y)q(w_k,w_k)^{-1}$ and continues the
movement controlled by the Markov chain $\eta$. All the newborn
particles behave themselves as independent copies of the parent
particle. We assume that the Markov chain $\eta$ is irreducible and
the matrix $Q$ is conservative, i.e. $\sum\nolimits_{y\in
S}{q(x,y)}=0$ as well as $q(x,y)\geq0$ for $x\neq y$ and
$q(x,x)\in(-\infty,0)$ for any $x\in S$. Denote by $f_k(s):={\sf
E}{s^{\xi_k}}$, $s\in[0,1]$, the probability generating function of
random variable $\xi_k$, $k=1,\ldots,N$. We will employ the standard
condition of the finite derivative $f_k'(1)$ existence, i.e. of
${\sf E}{\xi_k}$ finiteness for each $k=1,\ldots,N$. Let $\mu(t)$ be
the total particles number existing in CBP at time $t\geq0$ and the
local particles numbers $\mu(t;y)$ are the quantities of particles
located at separate points $y\in S$ at time $t$.

To formulate the main results of the paper let us introduce
additional notation. We temporarily forget that there are catalysts
at some elements of set $S$ and consider only the movement of a
particle on set $S$ in accordance with Markov chain $\eta$ with
generator $Q$ and starting point $x$. Let $_H\overline{\tau}_{x,y}$
be the time elapsed from the exit moment of this Markov chain out of
starting state $x$ till the first hitting point $y$ by the
particle whenever the particle trajectory does not pass the set
$H\subset S$. Otherwise, we put $_H\overline{\tau}_{x,y}=\infty$.
Extended random variable $_H\overline{\tau}_{x,y}$ is called
\emph{hitting time} of state $y$ \emph{under taboo} on set $H$ after
exit out of starting state $x$ (see, e.g., \cite{Chung}, Part~II, Section~11). Denote by $_H\overline{F}_{x,y}(t)$, $t\geq0$, the improper
cumulative distribution function of this extended random variable.
Mainly we will be interested in the situation when $H=W_k$ where
$W_k:=W\setminus\{w_k\}$, $k=1,\ldots,N$. Further
$$F^{\ast}(\lambda):=\int\nolimits_{0-}^{\infty}{e^{-\lambda t}\,d{F(t)}},\quad\lambda\geq0,$$ denotes
the Laplace transform of cumulative distribution function $F(t)$,
$t\geq0$, with support located on non-negative semi-axis. For
$j,k=1,\ldots,N$ and $t\geq0$ set $G_j(t):=1-e^{-\beta_j t}$,
$$G_{j,k}(t):=\beta_j\int\nolimits_0^{t}{{_{W_k}\overline{F}_{w_j,w_k}(t-u)}e^{-\beta_j u}\,du}$$
and
$${_{W_k}F_{x,w_k}(t)}:=-q(x,x)\int\nolimits_0^t{{_{W_k}\overline{F}_{x,w_k}(t-u)}e^{q(x,x)u}\,du}.$$
At last, let
$$_{W_k}\overline{F}_{x,w_k}(\infty)={_{W_k}F_{x,w_k}(\infty)}:=\lim_{t\to\infty}{_{W_k}F_{x,w_k}(t)}.$$

In \cite{B_TVP_14} there was introduced a matrix function
$D(\lambda)$ which is irreducible matrix of size $N\times N$ for
each $\lambda\geq0$. Namely,
$D(\lambda)=(d_{i,j}(\lambda))_{i,j=1}^N$ where
$$d_{i,j}(\lambda)=\delta_{i,j}\alpha_if_i'(1)G^{\ast}_i(\lambda)+(1-\alpha_i)G^{\ast}_i(\lambda)
{_{W_j}\overline{F}^{\ast}_{w_i,w_j}(\lambda)}$$ and $\delta_{i,j}$
is the Kronecker delta. According to definition~$1$ in
\cite{B_TVP_14} CBP is called supercritical if the Perron root (i.e.
positive eigenvalue being the spectral radius) $\rho(D(0))$ of the
matrix $D(0)$ is greater than $1$. Then in view of monotonicity of
all elements of matrix function $D(\cdot)$ there exists the solution
$\nu>0$ of equation $\rho(D(\lambda))=1$. As Theorem $1$ in
\cite{B_TVP_14} shows, just this positive number $\nu$ specifies the
rate of exponential growth of the mean total and local particles
numbers (in the literature devoted to population dynamics and
branching processes one traditionally speaks of Malthusian
parameter). More exactly, ${\sf E}_{x}\mu(t)\sim A(x)e^{\nu t}$ and
${\sf E}_x\mu(t;y)\sim a(x,y)e^{\nu t}$ as $t\to\infty$ where the
index $x\in S$ stands for the starting point of CBP. The explicit
formulae for functions $A(\cdot)$ and $a(\cdot,\cdot)$ are given in
\cite{B_TVP_14}.

Theorem \ref{q} is devoted to solution of the problem of finding the
global extinction probability $q(x):={\sf
P}_x(\lim_{t\to\infty}\mu(t)=0)=\lim_{t\to\infty}{\sf
P}_x(\mu(t)=0)$, $x\in S$, for the model of CBP.

\begin{Thm}\label{q}
For the global extinction probability $q(x)$ when $x\in S\setminus
W$ there exists a representation
\begin{equation}\label{eq.q}
q(x)=\sum_{k=1}^N{_{W_k}F_{x,w_k}(\infty)q(w_k)}
\end{equation}
where the values $q(w_j)$, $j=1,\ldots,N$, satisfy with the
following system of  equations
\begin{equation}\label{s_eq.q}
q(w_j)=\alpha_j f_j(q(w_j))+(1-\alpha_j)\sum_{k=1}^N{_{W_k}F_{w_j,w_k}(\infty)q(w_k)}.
\end{equation}
In addition, vector $(q(w_1),\ldots,q(w_N))$ is component-wise the
least root of equations system \eqref{s_eq.q} in the cube $[0,1]^N$.
Moreover, if the Markov chain $\eta$ is recurrent then $q(x)=1$,
$x\in S$, or $q(x)<1$, $x\in S$, whenever $\rho(D(0))\leq1$ or
$\rho(D(0))>1$, respectively. If the Markov chain $\eta$ is
transient then $q(x)<1$ for all $x\in S$.
\end{Thm}

One may oppose the global extinction probability $q(x)$ to the local
extinction probability $Q(x,y)={\sf
P}_x\left(\limsup_{t\to\infty}\mu(t;y)=0\right)$, for $x,y\in S$.
The following theorem shows that in fact the function $Q(x,y)$ does
not depend on variable $y$.

\begin{Thm}\label{Q}
Equality $Q(x,y)=Q(x)$ holds for any $y\in S$ where the function
$Q(x)$ when $x\in S\setminus W$ is of the form
\begin{equation}\label{eq.Q}
Q(x)=\sum_{k=1}^N{_{W_k}F_{x,w_k}(\infty)Q(w_k)}+1-\sum_{k=1}^N{_{W_k}F_{x,w_k}(\infty)}
\end{equation}
and the values $Q(w_j)$, $j=1,\ldots,N$, are the least solution to the equations system
\begin{eqnarray}\label{s_eq.Q}
Q(w_j)=\alpha_j f_j(Q(w_j))&+&(1-\alpha_j)\sum_{k=1}^N{_{W_k}F_{w_j,w_k}(\infty)Q(w_k)}\nonumber\\
&+&(1-\alpha_j)\left(1-\sum_{k=1}^N{_{W_k}F_{w_j,w_k}(\infty)}\right)
\end{eqnarray}
in the cube $[0,1]^N$. Moreover, if $\rho(D(0))\leq1$ then $Q(x)=1$ for all $x\in S$. Whenever $\rho(D(0))>1$
one has $Q(x)<1$ for each $x\in S$.
\end{Thm}

Applying Theorems 4 and 5 in \cite{Sev}, Chapter~5, Section~1, to the
auxiliary Bellman-Harris processes constructed in \cite{B_TVP_14} we
come to the statements of Theorems \ref{q} and \ref{Q}. Clearly
$0\leq q(x)\leq Q(x)\leq 1$, $x\in S$. In view of the explicit form
of relations \eqref{eq.q}-\eqref{s_eq.Q} we conclude that if the
strict inequality $\sum_{k=1}^N{_{W_k}F_{x,w_k}(\infty)}<1$ is
satisfied at least for some $x\in S$, i.e. the Markov chain $\eta$
is transient, then $q(x)<Q(x)$ for all $x\in S$. Otherwise, i.e. if
the Markov chain $\eta$ is recurrent, the relations for $q(\cdot)$
and $Q(\cdot)$ coincide and whence $q(x)=Q(x)$ for all $x$. In the
terms of paper \cite{Bertacchi_Zucca_14} the aforesaid means that,
for transient Markov chain $\eta$ and $\rho(D(0))\leq1$, we deal
with the pure global survival phase of CBP and in the case of
recurrent Markov chain $\eta$ and $\rho(D(0))>1$ the strong local
survival of CBP is observed.

Now we pass to considering the problem of the population growth rate
in the case of global and local survival, i.e. whenever
$\rho(D(0))>1$. Let ${\bf u}=(u_1,\ldots,u_N)$ be the right
eigenvector of matrix $D(\nu)$ corresponding to the Perron root
$\rho(D(\nu))$ equal to $1$ where $u_k>0$, $k=1,\ldots,N$, and
$\sum_{k=1}^Nu_k=1$. It should be noted that, by virtue of the
Perron-Frobenius theorem (see, e.g., \cite{Sev}, Chapter IV, Section~5),
such eigenvector can always be found since the matrix $D(\nu)$ is
irreducible according to Lemma 1 in \cite{B_TVP_14}. Recall the
definition of matrix function $D(x;\lambda)$, $x\notin W$,
$\lambda\geq0$, introduced in \cite{B_TVP_14} while proving Theorem
1, case 2. For $x\notin W$, set $w_{N+1}=x$, ${W(x):=W\cup\{x\}}$
and ${W_i(x):=W(x)\setminus\{w_i\}}$, ${i=1,\ldots,N+1}$. Then the
matrix ${D(x;\lambda)=(d_{i,j}(x;\lambda))_{i,j=1}^{N+1}}$ has
elements
$$d_{i,j}(x;\lambda):=\delta_{i,j}\,\alpha_i f_i'(1)G^\ast_i(\lambda)+(1-\alpha_i)G^\ast_i(\lambda){_{W_j(x)}\overline{F}^\ast_{w_i,w_j}(\lambda)},\quad\lambda\geq0.$$ Here $\alpha_{N+1}=0$, ${f'_{N+1}(1)=0}$ and
$G_{N+1}(t):=1-e^{q(x,x)t}$, ${t\geq0}$. By Lemma 3 in
\cite{B_TVP_14} the matrix $D(x;\lambda)$ is irreducible and
$\rho(D(x;\nu))=1$. Define ${\bf u}(x)=(u_1(x),\ldots,u_{N+1}(x))$
to be the right eigenvector of the matrix $D(x;\nu)$, corresponding
to the Perron root $\rho(D(x;\nu))$ equal to $1$, such that
$u_k(x)>0$, $k=1,\ldots,N+1$, and $\sum_{k=1}^{N}u_k(x)=1$. Take
notice that bearing on the proof of Lemma~$3$ in \cite{B_TVP_14} it
is not difficult to verify equalities $u_i(x)=u_i$ for all
$i=1,\ldots,N$ and $x\in S\setminus W$. Set $c(x):=u^{-1}_k$
whenever $x=w_k$ for some $k=1,\ldots,N$ and $c(x):=u^{-1}_{N+1}(x)$
whenever $x\in S\setminus W$. Let also symbols ${\bf 0}$ and ${\bf
1}$ denote vectors with zeros and units, respectively, as all
components, the dimension of vectors being contextually clear.

In the given paper, we handle three forms of random variables
convergence, viz., almost surely, in probability and in distribution
which are denoted by $\xrightarrow{\mbox{\it a.s.}}$,
$\xrightarrow{\sf P}$ and $\xrightarrow{d}$, respectively. Theorem
\ref{strong} describes strong convergence of vectors of the
normalized total and local particles numbers in CBP as time $t$
grows to infinity. The results on asymptotic behavior of the
normalizing means obtained in Theorem 1 in \cite{B_TVP_14} are
implicitly employed in the proofs of Theorems \ref{strong} and
\ref{weak}.

\begin{Thm}\label{strong}
Let a supercritical CBP start at point $x\in S$. Assume that the
elements of matrix $Q$ are uniformly bounded, i.e. for all
${z_1,z_2\in S}$ and some constant $C>0$ one has ${|q(z_1,z_2)|<C}$.
Let also ${\sf E}\xi^2_k<\infty$ for each $k=1,\ldots,N$. Then there
exists non-degenerate random variable $\zeta$ such that for any
$n\in\mathbb{N}$ and $y_1,\ldots,y_n\in S$ the following relation is
valid
\begin{equation}\label{T:strong}
\left(\frac{\mu(t)}{{\sf E}_x\mu(t)},\frac{\mu(t;y_1)}{{\sf E}_x\mu(t;y_1)},\ldots,
\frac{\mu(t;y_n)}{{\sf
E}_x\mu(t;y_n)}\right)\xrightarrow{\mbox{\it a.s.}}c(x)\zeta {\bf
1},\quad t\to\infty.
\end{equation}
\end{Thm}

If the Markov chain $\eta$ is recurrent, the assertion of
Theorem~\ref{strong} ensues from the construction of auxiliary
Bellman-Harris branching process in \cite{B_TVP_14} and applying to
it Theorem 4.1 in \cite{Mode_2_68} with due regard of  the fact that
under the conditions of Theorem \ref{strong} the distributions of
life-lengths of particles of different types in the auxiliary
Bellman-Harris process are absolutely continuous and have bounded
densities. The latter observation is true on account of the argument
in the final part of paper \cite{B_SPL_14}.

If the Markov chain $\eta$ is transient then one has to apply the
proof scheme of Theorem 4.1 in \cite{Mode_2_68} to the auxiliary
Bellman-Harris process with final type of particles which was
introduced in \cite{B_TVP_14} while establishing the case 1 (step 6)
of Theorem 1.

Now we relax the restrictions on the studied process. For
$z\in\mathbb{R}$, set $\log^{+}z:=\log(max\{z,1\})$.

\begin{Thm}\label{weak}
Let a supercritical CBP start at point $x\in S$. Then for each
$n\in\mathbb{N}$ and any $y_1,\ldots,y_n\in S$ the following
alternative is true.
\begin{enumerate}
  \item If ${\sf E}{\xi_k}\log^+{\xi_k}=\infty$ for some $k\in\{1,\ldots,N\}$ then
  \begin{equation}\label{T:weak_P}
  \left(\frac{\mu(t)}{{\sf E}_x\mu(t)},\frac{\mu(t;y_1)}{{\sf E}_x\mu(t;y_1)},\ldots,\frac{\mu(t;y_n)}{{\sf E}_x\mu(t;y_n)}\right)
  \xrightarrow{\mbox{\sf P}}{\bf 0},\quad t\to\infty.
  \end{equation}
  \item If ${\sf E}{\xi_k}\log^+{\xi_k}<\infty$ for all $k=1,\ldots,N$ then
  \begin{equation}\label{T:weak_d}
  \left(\frac{\mu(t)}{{\sf E}_x\mu(t)},\frac{\mu(t;y_1)}{{\sf E}_x\mu(t;y_1)},\ldots,\frac{\mu(t;y_n)}{{\sf E}_x\mu(t;y_n)}\right)
  \xrightarrow{\mbox{d}}c(x)\zeta{\bf 1},\quad t\to\infty.
  \end{equation}
\end{enumerate}
Here $\zeta$ is a non-degenerate random variable with the following properties.
\begin{description}
  \item[(i)] ${\sf E}_x\zeta=c(x)^{-1}$.
  \item[(ii)] ${\sf P}_x(\zeta=0)={\sf P}_x\left(\limsup_{t\to\infty}\mu(t;y)=0\right)$ for any $y\in S$.
  \item[(iii)] The Laplace transform $\varphi(\lambda;x):={\sf E}_xe^{-\lambda\zeta}$, $\lambda\geq0$, $x\in S$,
  of random variable $\zeta$ for $x\in S\setminus W$ is of the form
  $$\varphi(\lambda;x)=\sum_{k=1}^N{\int\nolimits_0^{\infty}{\varphi(\lambda e^{-\nu u};w_k)\,d {_{W_k}F_{x,w_k}(u)}}}+1-\sum_{k=1}^N{_{W_k}F_{x,w_k}(\infty)}$$
  where functions $\varphi(\cdot;w_j)$, $j=1,\ldots,N$, satisfy the system of integral equations
  \begin{eqnarray*}
  \varphi(\lambda;w_j)&=&\alpha_j\int\nolimits_0^{\infty}{f_j(\varphi(\lambda e^{-\nu u};w_j))\,dG_j(u)}\\
  &+&(1-\alpha_j)\sum_{k=1}^N{\int\nolimits_0^{\infty}{\varphi(\lambda e^{-\nu u};w_k)\,dG_{j,k}(u)}}\\
  &+&(1-\alpha_j)\left(1-\sum_{k=1}^N{_{W_k}\overline{F}_{w_j,w_k}(\infty)}\right).
  \end{eqnarray*}
  \item[(iv)] The conditional distribution of $\zeta$ under condition of CBP start at point $x\in S$ is
  absolutely continuous on the positive semi-axis and has continuous density function.
\end{description}
\end{Thm}

Note that in relation \eqref{T:weak_d} of Theorem \ref{weak} exactly
the convergence of vectors is essential whereas in formulae
\eqref{T:strong} and \eqref{T:weak_P} the vectors convergence is
tantamount to convergence of their components.

If the Markov chain $\eta$ is recurrent then Theorem \ref{weak} is
proven with the help of Theorem 1.1 in \cite{Kaplan_75}, applied to
the constructed in \cite{B_TVP_14} auxiliary Bellman-Harris process.
However, if the Markov chain $\eta$ is transient, the statement of
Theorem \ref{weak} is established by applying the proof scheme of
Theorem 1.1 in \cite{Kaplan_75} to the Bellman-Harris process with
final type of particles introduced while proving the case 1 (step 6)
of Theorem 1 in \cite{B_TVP_14}.

Thus, under weak conditions new asymptotic properties of CBP are investigated.

\end{document}